\RequirePackage{ifpdf}
\ifpdf 
\documentclass[pdftex]{sigma}
\else
\documentclass{sigma}
\fi

\begin{document}

\renewcommand{\PaperNumber}{016}

\FirstPageHeading

\ShortArticleName{Representations and Discrete $q$-Ultraspherical
Polynomials}

\ArticleName{Representations of the Quantum Algebra $\boldsymbol{{\rm su}_q(1,1)}$\\
and Discrete $\boldsymbol{q}$-Ultraspherical Polynomials}

\Author{Valentyna GROZA}

\AuthorNameForHeading{V. Groza}

\Address{National Aviation University, 1 Komarov Ave.,
Kyiv, 03058 Ukraine}

\Email{\href{mailto:groza@i.com.ua}{groza@i.com.ua}}

\ArticleDates{Received September 16, 2005, in final form November 09, 2005;
Published online November 15, 2005}

\Abstract{We derive orthogonality relations for discrete
$q$-ultraspherical polynomials and their duals by means of
operators of representations of the quantum algebra ${\rm
su}_q(1,1)$. Spectra and eigenfunctions of these operators are
found explicitly. These eigenfunctions, when normalized, form an
orthonormal basis in the representation space.}


\Classification{17B37; 33D45}

\section[Representations of ${\rm su}_q(1,1)$ with lowest
weights]{Representations of $\boldsymbol{{\rm su}_q(1,1)}$ with lowest
weights}

The aim of this paper is to study orthogonality relations for the
discrete $q$-ultraspherical polynomials and their duals by means
of operators of representations of the quantum algebra ${\rm
su}_q(1,1)$.

Throughout the sequel we always assume that $q$ is a fixed
positive number such that $q<1$. We use (without additional
explanation) notations of the theory of special functions and the
standard $q$-analysis (see, for example, \cite{[1]}).

The quantum algebra ${\rm su}_q(1,1)$ is defined as an associative
algebra, generated by the elements $J_+$, $J_-$, $q^{J_0}$ and
$q^{-J_0}$, subject to the defining relations
\begin{gather*}
q^{J_0}q^{-J_0}=q^{-J_0}q^{J_0}=1,
\\
q^{J_0}J_{\pm}q^{-J_0}=q^{\pm 1}J_{\pm},
\\
 [J_-,J_+] = {q^{J_0}-q^{-J_0}\over q^{1/2}-q^{-1/2}}  ,
\end{gather*}
and the involution relations $(q^{J_0})^* = q^{J_0}$ and $J_+^* =
J_-$. (We have replaced $J_-$ by $-J_-$ in the common definition
of the algebra $U_q({\rm sl}_{2})$; see \cite[Chapter 3]{[2]}.)

We are interested in representations of ${\rm su}_q(1,1)$ with
lowest weights. These irreducible representations are denoted by
$T^+_l$, where $l$ is a lowest weight, which can be {\it any
complex number} (see, for example, \cite{[3]}). They act on the
Hilbert space ${\cal H}$ with the orthonormal basis $|n\rangle$,
$n=0,1,2,\ldots$. The representation $T^+_{l}$ can be given in the
basis $|n\rangle$, $n=0,1,2,\ldots$, by the formulas
\begin{gather*}
q^{\pm J_0}\, |n\rangle  = q^{\pm (l + n)}\,|n\rangle ,
\\
J_+\,|n\rangle =\frac{q^{-(n+l-1/2)/2}}{1-q}
\sqrt{(1-q^{n+1})(1-q^{2l+n})} |n+1\rangle,
\\
J_-\, |n\rangle  = \frac{q^{-(n+l-3/2)/2}}{1-q}
\sqrt{(1-q^{n})(1-q^{2l+n-1})} |n-1\rangle .
\end{gather*}

For positive values of $l$ the representations $T^+_l$ are
$*$-representations. For studying discrete $q$-ultraspherical
polynomials we use the representations $T^+_l$ for which
$q^{2l-1}=-a$, $a>0$. They are not $*$-representations. But we
shall use operators of these representations which are symmetric
or self-adjoint. Note that $q^l$ is a pure imaginary number.

\section[Discrete $q$-ultraspherical polynomials and their
duals]{Discrete $\boldsymbol{q}$-ultraspherical polynomials and their
duals}

There are two types of discrete $q$-ultraspherical polynomials
\cite{[4]}. The first type, denoted as $C_n^{(a)}(x;q)$, $a>0$, is
a particular case of the well-known big $q$-Jacobi polynomials.
For this reason, we do not consider them in this paper. The second
type of discrete $q$-ultraspherical polynomials, denoted as
$\tilde C_n^{(a)}(x;q)$, $a>0$, is given by the formula
 \begin{gather}
 \tilde C_n^{(a)}(x;q)=(-{\rm i})^n C_n^{(-a)}({\rm i}x;q)= (-{\rm i})^n
{}_3\phi_2\left( q^{-n}, -aq^{n+1}, {\rm i}x;\,
   {\rm i}\sqrt{a}q, -{\rm i}\sqrt{a}q ;\, q,q \right) . \label{eq1}
 \end{gather}
 (Here and everywhere below under $\sqrt{a}$, $a>0$, we
understand a positive value of the root.)

The polynomials $\tilde C_n^{(a)}(x;q)$ satisfy the recurrence
relation
 \begin{gather}
 x \tilde C_n^{(a)}(x;q)= A_n \, \tilde C_{n+1}^{(a)}(x;q) +
C_n\, \tilde C_{n-1}^{(a)}(x;q),  \label{eq2}
 \end{gather}
where
\begin{gather*}
A_n=\frac{1+aq^{n+1}}{1+aq^{2n+1}} ,
\qquad
C_n=A_n-1=\frac{aq^{n+1}(1-q^n)}{1+aq^{2n+1}}.
\end{gather*}
Note that $ A_n \ge 1$ and, hence, coefficients in the recurrence
relation~\eqref{eq2} satisfy the conditions $ A_n C_{n+1}>0$ of
Favard's characterization theorem for $n=0,1,2,\ldots$. This means
that these polynomials are orthogonal with respect to a positive
measure. Orthogonality relation for them is derived in~\cite{[4]}.
We give here an approach to this orthogonality by means of
operators of representations $T^+_l$ of ${\rm su}_q(1,1)$.

Dual to the polynomials $C_n^{(a)}(x;q)$ are the polynomials
$D_n^{(a)}(\mu (x;a)|q)$, where $\mu(x;a)=q^{-x}+aq^{x+1}$. These
polynomials are a particular case of the dual big $q$-Jacobi
polynomials, studied in~\cite{[5]}, and we do not consider them.
Dual to the polynomials $\tilde C_n^{(a)}(x;q)$ are the
polynomials
 \begin{gather}
 \tilde D_n^{(a)}(\mu (x;-a)|q):= {}_3\phi_2\left(
q^{-x},-aq^{x+1},q^{-n};\,  {\rm i}\sqrt{a}q, -{\rm i}\sqrt{a}q;\,
 q,-q^{n+1} \right). \label{eq3}
 \end{gather}

For $a>0$ these polynomials satisfy the conditions of Favard's
theorem and, therefore, they are orthogonal. We derive an
orthogonality relation for them by means of operators of
representations $T^+_l$ of ${\rm su}_q(1,1)$.

\section[Representation operators $I$ and $J$]{Representation operators $\boldsymbol{I}$ and $\boldsymbol{J}$}

Let $T^+_l$ be the irreducible representation of ${\rm su}_q(1,1)$
with lowest weight $l$ such that $q^{2l-1}=-a$, $a>0$ (note that
$a$ can take any positive value). We consider the operator
 \begin{gather}
 I := \alpha\,q^{J_0/4}(J_+\, A +A\, J_-)\,q^{J_0/4}  \label{eq4}
 \end{gather}
of the representation $T^+_l$, where $\alpha =(a^2/q)^{1/2}(1-q)$
and
 \[
A=\frac{q^{(J_0-l+2)/2} \sqrt{(1-a^2q^{J_0-l+1})(1+aq^{J_0-l+1})}}
{\sqrt{(1+a^2q^{2J_0-2l+1})(1+a^2q^{2J_0-2l+2})(1+a^2q^{2J_0-2l+3})}}.
\]

We have the following formula for the symmetric operator $I$:
 \begin{gather}
I\, |n\rangle=a_n |n+1\rangle+a_{n-1}|n-1\rangle , \label{eq5}
 \\
 a_{n-1}=\big(a^2q^{n+1}\big)^{1/2}\, \left( \frac{
(1-q^{n})(1+a^2q^{n})}{ (1+a^2q^{2n-1})(1+a^2q^{2n+1})}
\right)^{1/2} .\nonumber
\end{gather}

The operator $I$ is bounded. We assume that it is defined on the
whole representation space~${\cal H}$. This means that $I$ is a
self-adjoint operator. Actually, $I$ is a Hilbert--Schmidt
operator since $a_{n+1}/a_n \to q^{1/2}$ when $n\to \infty$. Thus,
a spectrum of $I$ is simple (since it is representable by a~Jacobi
matrix with $a_n\ne 0$), discrete and have a single accumulation
point at~0 (see \cite[Chapter~VII]{[6]}).

To find eigenvectors $\psi_\lambda $ of the operator $I$, $I
\psi_\lambda =\lambda \psi_\lambda $, we set
 \[  
\psi_\lambda =\sum _{n} \beta_n(\lambda)|n\rangle.
 \]
Acting by $I$ upon both sides of this relation, one derives
 \[
\sum_n \beta_n(\lambda) (a_n |n+1\rangle+a_{n-1}|n-1\rangle )
=\lambda \sum \beta_n(\lambda)|n\rangle ,
 \]
where $a_n$ are the same as in~\eqref{eq5}. Collecting in this
identity factors at $|n\rangle$ with fixed $n$, we obtain the
recurrence relation for the coefficients $\beta_n(\lambda)$: $
a_n\beta_{n+1}(\lambda) +a_{n-1}\beta_{n-1}(\lambda) = \lambda
\beta_{n}(\lambda)$. Making the substitution
 \[
\beta_{n}(\lambda)=\left[\big(a^2q;q\big)_n \big(1+a^2q^{2n+1}\big)/(q;q)_n
\big(1+a^2q\big)a^{2n}\right]^{1/2}q^{-n(n+3)/4} \beta'_{n}(\lambda)
 \]
we reduce this relation to the following one
 \[
A_n \beta'_{n+1}(\lambda)+ C_n \beta'_{n-1}(\lambda) =\lambda
\beta'_{n}(\lambda) ,
 \]
where
 \[
 A_n=\frac{1+a^2q^{n+1}}{1+a^2q^{2n+1}},
  \qquad 
 C_n=\frac{a^2q^{n+1}(1-q^{n})}{1+a^2q^{2n+1}}.
 \]
 It is the recurrence
relation~\eqref{eq2} for the discrete $q$-ultraspherical
polynomials $\tilde C^{(a^2)}_n(\lambda;q)$. Therefore,
$\beta'_n(\lambda )=\tilde C^{(a^2)}_n(\lambda;q)$ and
 \begin{gather}
 \beta_n(\lambda )=\left(\frac{(-a^2q;q)_n\,(1+a^2q^{2n+1})}
{(q;q)_n\, (1+a^2q) a^{2n}} \right)^{1/2}q^{-n(n+3)/4}\, \tilde
C^{(a^2)}_n(\lambda;q). \label{eq6}
 \end{gather}
For the eigenfunctions $\psi _\lambda(x)$ we have the expansion
 \begin{gather}
 \psi _\lambda(x)=\sum_{n=0}^\infty
\left(\frac{(-a^2q;q)_n\,(1+a^2q^{2n+1})} {(q;q)_n\, (1+a^2q)
a^{2n}} \right)^{1/2}q^{-n(n+3)/4} \tilde C^{(a^2)}_n(\lambda;q)
 |n\rangle. \label{eq7}
 \end{gather}
Since a spectrum of the operator $I$ is discrete, only a discrete
set of these functions belongs to the Hilbert space ${\cal H}$ and
this discrete set determines the spectrum of $I$.

We intend to study the spectrum of $I$. It can be done by using
the operator
 \[
J:= q^{-J_0+l} - a^2\,q^{J_0-l+1}.
 \]
In order to determine how this operator acts upon the eigenvectors
$\psi _\lambda$, one can use the $q$-difference equation
 \begin{gather}
 \big(q^{-n}-a^2q^{n+1}\big)\, \tilde C^{(a^2)}_n(\lambda;q)= -a^2q
\lambda^{-2}\big( \lambda^2 +1\big)\,\tilde C^{(a^2)}_n(q\lambda;q)
\nonumber\\
\qquad{} + \lambda^{-2}a^2q(1+q)\, \tilde
C^{(a^2)}_n(\lambda;q)+ \lambda^{-2}\big(\lambda^2-a^2q^2\big)\tilde
C^{(a^2)}_n\big(q^{-1}\lambda;q\big)
 \label{eq8}
 \end{gather}
for the discrete $q$-polynomials polynomials. Multiply both sides
of~\eqref{eq8} by $d_n\,|n\rangle$, where $d_n$ are the
coefficients of $\tilde C^{(a^2)}_n(\lambda;q)$ in the expression~\eqref{eq6} 
for the coefficients $\beta_n(\lambda)$, and sum up
over $n$. Taking into account formula~\eqref{eq7} and the fact
that $J\,|n\rangle=(q^{-n}-a^2\,q^{n+1})\,|n\rangle$, one obtains
the relation
 \begin{gather}
 J\, \psi _{\lambda}=- a^2q\lambda^{-2}\big( \lambda^2 +1\big)
\psi_{q\lambda}  + \lambda^{-2}a^2q(1+q)\psi _{\lambda} +
\lambda^{-2}\big(\lambda^2-a^2q^2\big)\psi _{q^{-1}\lambda}  \label{eq9}
 \end{gather}
which is used below.

\section[Spectrum of $I$ and orthogonality of discrete
$q$-ultraspherical polynomials]{Spectrum of $\boldsymbol{I}$ and orthogonality of discrete\\
$\boldsymbol{q}$-ultraspherical polynomials}

Let us analyse a form of the spectrum of $I$ by using the
representations $T^+_l$ of the algebra ${\rm su}_q(1,1)$ and the
method of paper \cite{[7]}. If $\lambda$ is a spectral point of
$I$, then (as it is easy to see from~\eqref{eq9}) a successive
action by the operator $J$ upon the eigenvector $\psi_\lambda$
leads to the vectors~$ \psi_{q^m\lambda}$, $m=0,\pm 1, \pm
2,\ldots$. However, since $I$ is a Hilbert--Schmidt operator, not
all these points can belong to the spectrum of $I$, since
$q^{-m}\lambda \to\infty$ when $m\to \infty$ if $\lambda\ne 0$.
This means that the coefficient $\lambda^{-2}(\lambda^2 -a^2q^2)$
at $\psi _{q^{-1}\lambda}$ in~\eqref{eq9} must vanish for some
eigenvalue $\lambda$. There are two such values of $\lambda$:
$\lambda = aq$ and $\lambda = -aq$. Let us show that both of these
points are spectral points of~$I$. We have
 \[ \tilde
C^{(a^2)}_n(aq;q)={}_2\phi_1 \big(q^{-n}, a^2q^{n+1} ;\, -aq; \, q,q \big)
=a^2 q^{n(n+1)}.
 \]
Likewise,
 \[
 \tilde C^{(a^2)}_n(-aq;q)=a^2 q^{n(n+1)}.
 \]
Hence, for the scalar product $\langle \psi_{aq},
\psi_{aq}\rangle$ in ${\cal H}$ we have the expression
\begin{gather} \sum_{n=0}^\infty\,\frac{(1+a^2q^{2n+1})\,(-a^2q;q)_n}
{(1+a^2q)(q;q)_n\, a^{2n}q^{n(n+3)/2}} \tilde C^{(a^2)}_n(aq;q)^2=
\frac{(-a^2q^2,-q;q)_\infty}{(-a^2q^2;q^2)_\infty}\nonumber\\
\qquad {} =\left(\frac {(-a^2q^3;q^2)_\infty}{(q;q^2)_\infty}\right) <\infty. \label{eq10}
 \end{gather}
Similarly,
 \[
 \langle \psi_{-aq},\psi_{-aq} \rangle
=(-a^2q^2,-1;q)_\infty /(-a^2q^2;q^2)_\infty <\infty .
 \]
Thus, the values $\lambda=aq$ and $\lambda=-aq$ are the spectral
points of $I$.

Let us find other spectral points of the operator $I$. Setting
$\lambda = aq$ in~\eqref{eq9}, we see that the operator $J$
transforms $\psi _{aq}$ into a linear combination of the vectors
$\psi_{aq^2}$ and $\psi _{aq}$. We have to show that $\psi_{aq^2}$
belongs to the Hilbert space ${\cal H}$, that is, that
 \[
\langle \psi _{aq^2} ,\psi _{aq^2} \rangle = \sum_{n=0}^\infty\,
\frac{(1+a^2q^{2n+1})\,(-a^2q;q)_n} {(1+a^2q)(q;q)_n\,
a^{2n}}\,q^{-n(n+3)/2}\, \tilde C^{(a^2)}_n(aq^2;q)^2 <\infty .
 \]
It is made in the same way as in the case of big $q$-Jacobi
polynomials in paper \cite{[7]}. The above inequality shows that
$\psi _{aq^2}$ is an eigenvector of $I$ and the point $aq^2$
belongs to the spectrum of $I$. Setting $\lambda = aq^2$ in~\eqref{eq9} 
and acting similarly, one obtains that~$\psi _{aq^3}$
is an eigenvector of~$I$ and the point $aq^3$ belongs to the
spectrum of~$I$. Repeating this procedure, one sees that~$\psi
_{aq^n}$, $n=1,2,\ldots$, are eigenvectors of~$I$ and the set
$aq^n$, $n=1,2,\ldots$, belongs to the spectrum of~$I$. Likewise,
one concludes that $\psi _{-aq^n}$, $n=1,2,\ldots$, are
eigenvectors of~$I$ and the set $-aq^n$, $n=1,2,\ldots$, belongs
to the spectrum of~$I$. Let us show that the operator $I$ has no
other spectral points.

The vectors $\psi _{aq^n}$ and $\psi _{-aq^n}$, $n=1,2,\ldots$,
are linearly independent elements of the representation space
${\cal H}$. Suppose that $aq^n$ and $-aq^n$, $n=1,2,\ldots$,
constitute the whole spectrum of $I$. Then the set of vectors
$\psi _{aq^n}$ and $\psi_{-aq^n}$, $n=1,2,\ldots$, is a basis of
${\cal H}$. Introducing the notations $\Xi _n:=\psi_{aq^{n+1}}$ and
$\Xi'_n:=\psi _{-aq^{n+1}}$, $n=0,1,2,\ldots$, we find from~\eqref{eq9} that
\begin{gather*}
J\, \Xi _n = -q^{-2n-1}\big(1+a^2q^{2(n+1)}\big) \Xi _{n+1} + d_n\,\Xi _n
- q^{-2n}\big(1-q^{2n}\big)\, \Xi _{n-1} ,
\\
J \,\Xi'_n = -q^{-2n-1}\big(1+a^2q^{2(n+1)}\big)\, \Xi' _{n+1} + d_n\,\Xi'
_n - q^{-2n}\big(1-q^{2n}\big)\, \Xi'_{n-1} ,
\end{gather*}
where $ d_n=  q^{-2n-1}(1+q)$.

As we see, the matrix of the operator $J$ in the basis $\Xi _n $,
$\Xi'_n$, $n=0,1,2,\ldots$, is not symmetric, although in the
initial basis $|n\rangle$, $n=0,1,2,\ldots$, it was symmetric. The
reason is that the matrix $M:=(a_{mn}\; a'_{m'n'})$ with entries
 \[
 a_{mn}:=\beta_m\big(aq^{n+1}\big),\qquad  a'_{m'n'}:=\beta_{m'}\big(-aq^{n'+1}\big),
\qquad m,n,m',n'=0,1,2,\ldots,
 \]
where $\beta_m(dq^{n+1})$, $d=\pm a$, are coefficients~\eqref{eq6} in the expansion
 \[
 \psi _{dq^{n+1}}=\sum _m \beta_m\big(dq^{n+1}\big)\,|n\rangle,
 \]
is not unitary. (This matrix $M$ is formed by adding the columns
of the matrix $(a'_{m'n'})$ to the columns of the matrix
$(a_{mn})$ from the right.) It maps the basis $\{ |n\rangle \}$
into the basis $\{\psi_{aq^{n+1}}, \psi _{-aq^{n+1}} \}$ in the
representation space. The nonunitarity of the matrix $M$ is
equivalent to the statement that the basis $\Xi _n$, $\Xi _n$,
$n=0,1,2,\ldots$, is not normalized. In order to normalize it we
have to multiply $\Xi _n$ by appropriate numbers $c_n$ and
$\Xi'_n$ by numbers $c'_n$. Let
 \[
 \hat\Xi _n = c_n\Xi _n,\qquad \hat\Xi'_n =c'_n\Xi_n,\qquad
n=0,1,2,\ldots ,
 \]
be a normalized basis. Then the operator $J$ is symmetric in this
basis and has the form
\begin{gather*}
J\,\hat\Xi _n = -c_{n+1}^{-1}c_nq^{-2n-1}\big(1+a^2q^{2(n+1)}\big)
\,\hat\Xi_{n+1} + d_n\, \hat\Xi _n - c_{n-1}^{-1}c_n
q^{-2n}\big(1-q^{2n}\big)\,\hat\Xi _{n-1},
\\
J\,\hat\Xi'_n = -{c'}_{n+1}^{-1}{c'}_nq^{-2n-1}\big(1+a^2q^{2(n+1)}\big)
\,\hat\Xi' _{n+1} + d_n \,\hat\Xi' _n  -{c'}_{n-1}^{-1}{c'}_n
q^{-2n} \big(1-q^{2n}\big)\,\hat\Xi'_{n-1} .
\end{gather*}
The symmetricity of the matrix of the operator $J$ in the basis
$\{ \hat\Xi _n,\hat\Xi'_n \}$ means that for coefficients~$c_n$ we
have the relation
 \[
 c_{n+1}^{-1}c_n q^{-2n-1}\big(1+a^2q^{2(n+1)}\big)
=c_{n}^{-1}c_{n+1} q^{-2n-2} \big(1-q^{2(n+1)}\big).
 \]
The relation for $c'_n$ coincides with this relation. Thus,
 \[
\frac{c_{n}}{c_{n-1}}=\frac{c'_{n}}{c'_{n-1}}=
\sqrt{\frac{q(1+a^2q^{2n})}{1-q^{2n}}}.
 \]
This means that
 \[
c_n= C\left(
\frac{q^{n}(-a^2q^2;q^2)_n}{(q^2;q^2)_n}\right)^{1/2},
 \qquad
c'_n= C'\left(
\frac{q^{n}(-a^2q^2;q^2)_n}{(q^2;q^2)_n}\right)^{1/2},
 \]
where $C$ and $C'$ are some constants.
\newpage

Therefore, in the expansions
 \[ \hat\Xi _n\equiv \sum _m \,{\hat
a}_{mn}\, |m\rangle ,\qquad \hat\Xi _n(x) \equiv\sum _m \,{\hat
a}'_{mn}\, |m\rangle
 \]
the matrix ${\hat M}:=({\hat a}_{mn}\ {\hat a}'_{mn})$ with
entries ${\hat a}_{mn}= c_n\,\beta _m(aq^n)$ and ${\hat a}'_{mn}=
c_n\,\beta _m(cq^n)$ is unitary, provided that the constants $C$
and $C'$ are appropriately chosen. In order to calculate these
constants, one can use the relations $\sum\limits_{m=0}^\infty |{\hat
a}_{mn}|^2=1$ and $\sum \limits_{m=0}^\infty |{\hat a}'_{mn}|^2=1$ for
$n=0$. Then these sums are multiples of the sum in~\eqref{eq10},
so we find that
 \[
C=C'=\left(\frac{(-a^2q^2;q^2)_\infty}{(-a^2q^2,-q;q)_\infty}\right)^{1/2}
=\left( \frac{(q;q^2)_\infty}{(-a^2q^3;q^2)_\infty}\right)^{1/2} .
 \]

The orthogonality of the matrix ${\hat M}\equiv ({\hat a}_{mn}\
{\hat a}'_{mn})$ means that
 \begin{gather}
 \sum _m {\hat a}_{mn}{\hat a}_{mn'}=\delta_{nn'},\qquad 
  \sum _m {\hat a}'_{mn}{\hat a}'_{mn'}=\delta_{nn'},\qquad \sum _m
{\hat a}_{mn}{\hat a}'_{mn'}=0,  \label{eq11}
 \\
 \sum _n ({\hat a}_{mn}{\hat a}_{m'n}+ {\hat
a}'_{mn} {\hat a}'_{m'n} ) =\delta_{mm'} .   \label{eq12}
 \end{gather}
Substituting the expressions for ${\hat a}_{mn}$ and ${\hat
a}'_{mn}$ into~\eqref{eq12}, one obtains the relation
  \begin{gather}
\sum_{n=0}^\infty\, \frac{(-a^2q^2;q^2)_n q^n}{(q^2;q^2)_n}\left[\tilde
C^{(a^2)}_m\big(aq^{n+1}\big) \tilde C^{(a^2)}_{m'}\big(aq^{n+1}\big)  +\tilde
C^{(a^2)}_m\big(-aq^{n+1}\big) \tilde C^{(a^2)}_{m'}\big(-aq^{n+1}\big)\right]
\nonumber\\
\qquad{}= \frac{(-a^2q^3;q^2)_\infty}{(q;q^2)_\infty}
\frac{(1+a^2q)(q;q)_ma^{2m}
q^{m(m+3)/2}}{(1+a^2q^{2m+1})(-a^2q;q)_m} \delta_{mm'} .
\label{eq13}
 \end{gather}
 This identity must give an orthogonality relation for the
discrete $q$-ultraspherical polynomials $\tilde
C^{(a^2)}_{m}(y)\equiv \tilde C^{(a^2)}_{m}(y;q)$. An only gap,
which appears here, is the following. We have assumed that the
points $aq^{n+1}$ and $-aq^{n+1}$, $n=0,1,2,\ldots$, exhaust the
whole spectrum of the operator~$I$. If the operator $I$ would have
other spectral points $x_k$, then on the left-hand side 
of~\eqref{eq13} would appear other summands $\mu_{x_k} \tilde
C^{(a^2)}_m(x_k;q) \tilde C^{(a^2)}_{m'}(x_k;q)$, which correspond
to these additional points. Let us show that these additional
summands do not appear. For this we set $m=m'=0$ in the 
relation~\eqref{eq13} with the additional summands. This results in the
equality
 \begin{gather}
 \frac{(-a^2q^2
;q^2)_\infty}{(-a^2q^2,-q;q)_\infty}
\sum_{n=0}^\infty\,
\frac{(aq,-aq;q)_n q^n}{(q,-q;q)_n} \nonumber\\
\qquad {} +\frac{(-a^2q^2 ;q)_\infty}{(-a^2q^2,-q;q)_\infty} \sum_{n=0}^\infty
\,\frac{(aq,-aq;q)_n q^n}{(q,-q;q)_n}
+ \sum_k\,\mu_{x_k} =1. \label{eq14}
 \end{gather}
 In order to show that $\sum\limits_k \mu_{x_k} = 0$, take into account
formula (2.10.13) in \cite{[1]}. By means of this formula it is
easy to show that the relation~\eqref{eq14} without the summand
$\sum\limits_k \mu_{x_k}$ is true. Therefore, in~\eqref{eq14} the sum
$\sum\limits_k \mu_{x_k}$ does really vanish and formula~\eqref{eq13}
gives an orthogonality relation for the discrete
$q$-ultraspherical polynomials.

The relation~\eqref{eq13} and the results of Chapter VII in
\cite{[6]} shows that {\it the spectrum of the operator~$I$
coincides with the set of points} $aq^{n+1}$, $-aq^{n+1}$,
$n=0,1,2,\ldots$.

\section[Dual discrete $q$-ultraspherical polynomials]{Dual discrete $\boldsymbol{q}$-ultraspherical polynomials}

Now we use the relations~\eqref{eq11}. They give the
orthogonality relation for the set of matrix ele\-ments~${\hat
a}_{mn}$ and ${\hat a}'_{mn}$, considered as functions of $m$. Up
to multiplicative factors, they coincide with the functions
 \begin{gather*}
F_n(x;a^2):={}_3\phi_2 \big(x,a^2q/x, aq^{n+1};\, {\rm i}aq, -{\rm
i}aq;\, q,q\big),\qquad n=0,1,2,\ldots,
 \\
F'_n(x;a^2):={}_3\phi_2 \big(x,a^2q/x, -aq^{n+1};\, {\rm i}aq, -{\rm
i}aq;\, q,q\big),\qquad n=0,1,2,\ldots ,
 \end{gather*}
considered on the corresponding sets of points.

Applying the relation (III.12) of Appendix III in \cite{[1]} we
express these functions in terms of dual $q$-ultraspherical
polynomials. Thus we obtain expressions for ${\hat a}_{mn}$ and
${\hat a}'_{mn}$ in terms of these polynomials. Substituting these
expressions into the relations~\eqref{eq11} we obtain the
following orthogonality relation for the polynomials $\tilde
D^{(a^2)}_n(\mu(m;a^2)|q)$:
\begin{gather}
\sum_{m=0}^\infty
\frac{(1+a^2q^{2m+1})(-a^2q;q)_{m}}{(1+a^2q)(q;q)_{m}}
q^{m(m-1)/2} \tilde D_n^{(a^2)}\big(\mu \big(m;-a^2\big)|q\big)\tilde
D_{n'}^{(a^2)}\big(\mu \big(m;-a^2\big)|q\big)
 \nonumber\\
\qquad {}=\frac{2(-a^2q^3;q^2)_\infty}{(q;q^2)_\infty}
 \frac{(q^2;q^2)_nq^{-n}}{(-a^2q^2;q^2)_n} \delta_{nn'} . \label{eq15}
 \end{gather}
 This orthogonality relation coincides with the sum of two
orthogonality relations~\eqref{eq9} and~\eqref{eq10} in~\cite{[4]}. 
The orthogonality measure in~\eqref{eq15} is not
extremal since it is a sum of two extremal measures.

In order to obtain orthogonality relations~\eqref{eq9} and~\eqref{eq10} 
of \cite{[4]}, from the very beginning, instead of
operator $I$, we have to consider operators $I_1$ and $I_2$ of the
representation $T^+_l$ of the algebra $U_{q^2}({\rm su}_{1,1})$,
which are appropriate for obtaining the orthogonality relations
for the sets of polynomials $\tilde
C^{(a)}_{2k}(\sqrt{a}q^{s+1};q)$, $k=0,1,2,\ldots$, and $\tilde
C^{(a)}_{2k+1}(\sqrt{a}q^{s+1};q)$, $k=0,1,2,\ldots$, from Section~2 in~\cite{[4]}. 
Then going to the orthogonality relations for
dual sets of polynomials in the same way as above, we obtain the
extremal orthogonality relations~\eqref{eq9} and~\eqref{eq10} of~\cite{[4]}.

\LastPageEnding

\end{document}